\documentclass[10pt, 14paper,reqno]{amsart}
\usepackage{amsmath,amsfonts,amssymb}
\usepackage{longtable}
\usepackage{color}
\theoremstyle{plain}
\newtheorem{theorem}{Theorem}

\theoremstyle{proof}
\theoremstyle{definition}

\theoremstyle{remark}

\theoremstyle{lamma}

\numberwithin{equation}{section}
\numberwithin{lemma}{section}
\numberwithin{theorem}{section}
\numberwithin{remark}{section}
\numberwithin{prop}{section}
\numberwithin{corollary}{section}

\usepackage{amsmath}
\usepackage{amsfonts}   
\usepackage{amssymb}
\usepackage{amssymb, amsmath, amsthm}
\usepackage[breaklinks]{hyperref}

\theoremstyle{thmrm}

\numberwithin{conjecture}{section}

\begin{document}
\title[A survey report]{\large{Jacobi sums and cyclotomic numbers: A survey report}}
\author{Md Helal Ahmed and Jagmohan Tanti}
\address{Md Helal Ahmed @ Department of Mathematics, Central University of Jharkhand, Ranchi-835205, India}
\email{ahmed.helal@cuj.ac.in}

\address{Jagmohan Tanti @ Department of Mathematics, Central University of Jharkhand, Ranchi-835205, India}
\email{jagmohan.t@gmail.com}.

\maketitle
\begin{abstract}
The determination of Jacobi sums, their congruences and cyclotomic numbers have been the object of attention for many years and there are large 
number of interesting results related to these in the literature. This survey aims at reviewing results concerning the diophantine systems for 
finding the cyclotomic numbers and coefficients of Jacobi sums and to indicate the current status of the problem. 
\end{abstract}

\section{Introduction}
The Jacobi's tremendous mathematical legacy has many contributions to the field of mathematics, among which are the Jacobi symbol, the Jacobi triple 
product, the Jacobian in the change of Variables theorem and the Jacobi elliptic functions. Among his multiple discoveries, Jacobi sums appear as one of
the most important findings. In any given finite field $\mathbb{F}_q$, Jacobi sums of order $e$ mainly depend on two parameters. Therefore, these values 
could be naturally assembled into a matrix of order $e$. Jacobi initially proposed these sums as mathematical objects, and for more certainty, he mailed 
them to Gauss in 1827 (see \cite{Berndt1,Jacobi1}). After ten years, Jacobi \cite{Jacobi2} published his findings including all the extensions provided 
by other scholars such as Cauchy, Gauss and Eisenstein. It is worth mentioning that while Gauss sums suffice for a proof of quadratic reciprocity, a 
demonstration of cubic reciprocity law along similar lines requires a foray into the realm of Jacobi sums. In order to prove biquadratic reciprocity,
Eisenstein \cite{Eisenstein1} formulated a generalization of Jacobi sums. As illustrated in \cite{Ireland1}, Jacobi sums could be used for estimating 
the number of integral solutions to congruences such as $x^{3}+y^{3}\equiv 1 \pmod p$. These estimates played a key role in the development of 
Weil conjectures \cite{Weil1}. Jacobi sums could be used for the determination of a number of solutions of diagonal equations over finite fields. 
Jacobi sums were also utilized in primality test by Adleman, Pomerance and Rumely \cite{Adleman1}. The Problem of congruences of Jacobi sums of order $e$ 
concerns to establish certain congruence conditions modulo an appropriate power of  
 $(1-\zeta_e)$ which is useful to determine an element (which is co-prime to $e$) in $\mathbb{Z}[\zeta_e]$  as a Jacobi sum of order $e$. 
 It is worth mentioning that congruence conditions for Jacobi sums play major role for determination of algebraic characterization/ diophantine systems 
 of Jacobi sums, hence of all Jacobi sums together with the absolute value and prime ideal decomposition of Jacobi sums.

Cyclotomic numbers are one of the most important objects in number theory and in other branches of mathematics. These number have been extensively used in
coding theory, cryptography and in other branches of information theory. One of the central problems in the study of these numbers is the determination of
all cyclotomic numbers of a specific order for a given field in terms of solutions of certain Diophantine system.
This problem has been treated by many mathematicians including C. F. Gauss who had determined all the cyclotomic numbers of order $3$ in the field 
$\mathbb{F}_q$ with prime $q\equiv 1\pmod 3$. Complete solutions to this cyclotomic number problem have been computed for some specific orders. 
For instance, the cyclotomic numbers of prime order $e$ 
in the finite field $\mathbb{F}_q$ with $q=p^r$ and $p\equiv 1 \pmod e$ have been investigated by many authors (see, \cite{Katre3} and the references 
therein). Cyclotomic numbers of order $e$ over the field $\mathbb{F}_q$ with characteristic $p$, in general, cannot 
be determined only in terms of $p$ and $e$, but that one requires a quadratic partition of $q$ too.

In this survey article, we discuss some interesting results concerning the Jacobi sums and its congruences, and cyclotomic numbers as well as the 
current status of the problem. Starting from Gauss, this topic has been studied extensively by many authors and thus there exist a large number of 
research articles. Due to the versatility, this survey may miss out some interesting references and thus some interesting results too and thus this 
article is never claimed to be a complete one. 

\section{Definitions and notations}
Let $e\geq2$ be an integer, $p$ a rational prime,  $q=p^{r}, r \in\mathbb{Z}^{+}$ and $q\equiv 1 \pmod{e}$. Let $\mathbb{F}_{q}$ 
be a finite field of $q$ elements. We can write $q=p^{r}=ek+1$ for some $k\in\mathbb{Z}^{+}$. Let $\gamma$ be a generator of 
the cyclic group $\mathbb{F}^{*}_{q}$ and $\zeta_e=exp(2\pi i/e)$. Also for $a\in\mathbb{F}_q^*$, ${\rm ind}_{\gamma}(a)$ is defined to be a positive integer 
$m\leq q-1$ such that $a=\gamma^m$.
Define a multiplicative character $\chi_e : \ \mathbb{F}^{*}_{q} \longrightarrow \mathbb{Q}(\zeta_e)$ by $\chi_e(\gamma)=\zeta_e$ and extend it on 
$\mathbb{F}_q$ 
by putting $\chi_e(0)=0$. 
For  integers $\displaystyle 0\leq i,j\leq e-1$, the Jacobi sum $J_{e}(i,j)$ is define by

$$J_{e}(i,j)= \sum_{v\in \mathbb{F}_{q}} \chi_e^{i}(v) \chi_e^{j}(v+1).$$

However in the literature a variation of Jacobi sums are also considered and is defined by 
$$J_e(\chi_e^i,\chi_e^j)=\sum_{v\in\mathbb{F}_q}\chi_e^i(v)\chi_e^j(1-v).$$\\ Observe that $J_e(i,j)=\chi_e^i(-1)J_e(\chi_e^i,\chi_e^j)$. When $q=2^{r}$, $\chi_e^i(-1)=\chi_e^i(1)=1$ and both the Jacobi sums coincide. Otherwise $\chi_e^i(-1)=(-1)^{ik}$ and hence the two Jacobi sums differ at most in sign.

For $0\leq a, b \leq e-1$, the cyclotomic number $(a,b)_{e}$ of order $e$ is defined as the number of solutions $(s,t)$ of the following:
\begin{equation} \label{1}
\gamma^{es+a}+\gamma^{et+b}+1\equiv 0 \pmod q; \ \ \ \ \ \ 0\leq s,t \leq k-1.
\end{equation} 
\begin{center}
or
\end{center}
One can define, for $0\leq a, b\leq e-1$, the cyclotomic numbers $(a, b)_e$ of order $e$ is as follows:
\begin{align*}
(a,b)_e:& =\#\{v\in\mathbb{F}_q|\chi_e(v)=\zeta_e^a,\,\,\chi_e(v+1)=\zeta_e^b\} \\ & =
\#\{v\in\mathbb{F}_q\setminus \{0,-1\}\mid {\rm ind}_{\gamma}v\equiv a \pmod e ,  {\ \rm ind}_{\gamma}(v+1)\equiv b \pmod e\}.
\end{align*}
The cyclotomic numbers $(a,b)_e$ and the Jacobi sums $J_e(i,j)$ are well connected  by the following relations \cite{Berndt1,Shirolkar1}: 
\begin{equation} \label{01}
 \sum_a\sum_b(a,b)_e\zeta_e^{ai+bj}=J_e(i,j),
 \end{equation} 
 and
 \begin{equation} \label{00}
 \sum_i\sum_j\zeta_e^{-(ai+bj)}J_e(i,j)=e^{2}(a,b)_e.
 \end{equation} 
\eqref{01} and \eqref{00} show that if we want to calculate all the cyclotomic numbers $(a,b)_e$ of order $e$, it is sufficient to calculate all the 
Jacobi sums $J_e(i,j)$ of the same order, and vice-versa.
\section{Properties of Jacobi sums and cyclotomic number} 
In the following theorem we state some standard results about Jacobi sums.
\begin{theorem} \cite{Acharya1,Helal1,Berndt1} (Elementary properties of Jacobi sums)
\\
\textbf{i}.
If $m+n+s\equiv 0 \ ($mod$\ e)$ then 
\\
$J_{e}(m,n)=J_{e}(s,n)=\chi_e^{s}(-1)J_{e}(s,m)=\chi_e^{s}(-1)J_{e}(n,m)=\chi_e^{m}(-1)J_{e}(m,s)$ $=\chi_e^{m}(-1)J_{e}(n,s).$
 \\ 	
In particular,
\begin{equation*}
J_{e}(1,m)=\chi_e(-1)J_{e}(1,s)=\chi_e(-1)J_{e}(1,e-m-1).
\end{equation*} 
\textbf{ii}. 
\ \ \ \ $J_{e}(0,j)= \begin{cases}
-1 \ \ \ \ \ if  \ j\not\equiv 0 \ ($mod$\ e) ,\\
q-2 \ \ if  \ j\equiv 0 \ ($mod$\ e).
\end{cases}$

\ \ \ \ $J_{e}(i,0)=-\chi_e^{i}(-1) \ if \  i\not\equiv 0 \ ($mod$\ e).$
\\ \\ 
\textbf{iii}.  \ \ \ \ Let $m+n\equiv 0 \ ($mod$\ e)$ \ but not both $m$ and $n$ zero $\pmod{e}$. Then $J_{e}(m,n)=-1.$
\\ \\ \textbf{iv}.  \ \  For $(k,e)=1$ and $\sigma_k$ a $\mathbb{Q}$ automorphism of $\mathbb{Q}(\zeta_e)$ with $\sigma_k(\zeta_e)=\zeta_e^k$, 
we have $\sigma_{k}J_{e}(m,n)=J_{e}(mk,nk)$. In particular, if $(m,e)=1,\ m^{-1}$ \ denotes the inverse of $m\pmod{e}$ then 
$\sigma_{m^{-1}}J_{e}(m,n)=J_{e}(1,nm^{-1})$.
\\ \\ \textbf{v}.  \ \ \ \ 
$J_{2e}(2m,2s)=J_{e}(m,n)$. 
\\  \textbf{vi}. \ \ \ \ $J_{e}(1,n) \overline{J_{e}(1,n)}= \begin{cases}
q \ \ \ \ \ if \ n\not\equiv 0, -1 \ ($mod$\ e),\\
1 \ \ \ \ \ if \ n\equiv 0, -1 \ ($mod$\ e).
\end{cases}$
\\ \\
\textbf{vii}. \ \ \ Let $m$, $n$, $s$ be integers and $l$ be an odd prime, such that $m+n \not\equiv 0 \ ($mod$\ 2l)$ and $m+s \not\equiv 0 \ ($mod$\ 2l)$. Then 
\begin{equation*}
J_{2l}(m,n) J_{2l}(m+n ,s) = \chi^{m}(-1)J_{2l}(m,s)J_{2l}(n,s+m).
\end{equation*}
\textbf{viii}. \ \ \ Let $m$, $n$, $s$ be integers  and $l$ be an odd prime, such that $m+n \not\equiv 0 \ ($mod$\ 2l^{2})$ and $m+s \not\equiv 0 \ ($mod$\ 2l^{2})$. Then 
\begin{equation*}
J_{2l^{2}}(m,n) J_{2l^{2}}(m+n ,s) = \chi^{m}(-1)J_{2l^{2}}(m,s)J_{2l^{2}}(n,s+m).
\end{equation*}
\end{theorem}
In the next theorem, we state some basic properties of the cyclotomic numbers of order $e$.
\begin{theorem} \cite{Berndt1}
The cyclotomic numbers of order $e$ have the have the following properties:\\
\textbf{i}. \ \ \ $(a,b)_{e}=(a' ,b')_{e}$ if $a\equiv a'\pmod{e}$ and $b\equiv b'\pmod{e}$. \\
\textbf{ii}. \ \ \ $(a,b)_{e}= (e-a,b-a)_{e}$ along with the following:
\begin{equation*}\label{2.1}
(a,b)_{e}=\begin{cases}
 (b,a)_{e}\hspace*{1.751cm} \text{ if }  k \text{ is  even  or  } q=2^r,\\
(b+\frac{e}{2},a+\frac{e}{2})_{e}\hspace*{3 mm} \text{ otherwise }.
\end{cases}. 
\end{equation*}
\\
\textbf{iii}. \ \ \ \begin{equation*} \label{2.2}
\sum_{a=0}^{e-1}\sum_{b=0}^{e-1}(a,b)_{e}=q-2,
\end{equation*} 
\\
\textbf{iv}. \ \ \ \begin{equation*} \label{2.3}
\sum_{b=0}^{e-1}(a,b)_{e}=k-n_{a},
\end{equation*}
where $n_{a}$ is given by 
\begin{equation*}
n_{a}=\begin{cases}
1 \quad \text{ if } a=0, 2\mid k \text{ or if } a=\frac{e}{2}, 2\nmid k;\\
0 \quad \text{ otherwise }.
\end{cases}
\end{equation*}
\\
\textbf{v}. \ \ \ \begin{equation} \label{2.4}
\sum_{a=0}^{e-1}(a,b)_{e}=\begin{cases}
k-1 \hspace*{3mm}\text{ if } b =0; \\ 
k \hspace*{1cm}\text{ if } 1\leq b \leq e-1.
\end{cases}
\end{equation} 
\\
\textbf{vi}. \ \ \ $(a,b)_{e}^{\prime}=$  $(r^{_{0}}a,r^{_{0}}b)_{e}$, \\
where the prime $(\prime)$ indicates that the cyclotomic number is taken with respect to the generator $\gamma^{r_{0}}$ in place of $\gamma$ in $\mathbb{F}_{q}^{*}$.
\end{theorem}

\section{Jacobi sums and it congruences}
Gauss theories represented the cornerstone of Jacobi sum findings. Many research work have been conducted by a number of mathematicians in an attempt 
to find out the Diophantine systems that characterize the coefficients of Jacobi sums, i.e. giving a Diophantine system whose unique solution provides 
the coefficients of a particular Jacobi sum. Jacobi sums are particularly used for obtaining the cyclotomic numbers of the same order and vice-versa 
(i.e. the cyclotomic numbers of order $e$ are known if one knows all the Jacobi sums of order $e$ and vice-versa). Evaluating all the Jacobi sums of
order $e$ is relatively intricate. A number of authors devoted for the evaluation of Jacobi sums with certain order. Obtaining the concerned relations 
helps in reducing the complexity of evaluating all Jacobi sums as well as the cyclotomic numbers. The evaluations and relationships of Jacobi sums of 
order $3, \ 4$ and $7$ were introduced by Jacobi himself in a letter \cite{Jacobi1} to Gauss in 1827.  Relationships between the sum of order $e$ for 
$e\leq 6$, $e=8,\ 10$ and $12$ were established by Dickson \cite{Dickson1}. In later stages, Muskat \cite{Muskat1} established the relation of order 
$12$ in terms of the fourth root of unity to resolve the sign of ambiguity. Dickson \cite{Dickson3} found specific relationships for sums of order
$15,\ 16, \ 20$ and $24$. Muskat \cite{Muskat1} developed Dickson’s work for $e = 15$ and $24$ and extended it to sums of order $30$. Complete methods of 
$e = 16$ and $20$ exist in Whiteman \cite{Whiteman5} and Muskat \cite{Muskat3} respectively. In fact, before Dickson's work 
\cite{Dickson1,Dickson2,Dickson3}, Western \cite{Western1} determined Jacobi sums of order $8,\ 9,$ and $16$. An important issue that should be borne 
in mind is that all theories showed that Jacobi sums of higher orders can be expressed in terms of Jacobi sums of lower orders. Dickson \cite{Dickson2} 
gave some particular relationships for sums of order $14$ and $22$. Muskat \cite{Muskat2} provided complete results for order $14$. Dickson 
\cite{Dickson3} also investigated sums of order $9$ and $18$, while Baumert and Fredrickson \cite{Baumert1} gave corrections to some of his results and 
removed the sign of ambiguity. Zee also found relationships for sums of order $13$ and $60$ in \cite{Zee1}, and investigated the sums of order $22$ 
in \cite{Zee2}. Relationships for order $21,\ 28,\ 39,\ 55$ and $56$ are provided in one of Muskat and Zee research works \cite{Muskat4}. 
Berndt and Evans \cite{Berndt2} obtained sums of orders $3,\ 4, \ 6,\ 8,\ 12,\ 20$ and $24$, and they also determined sums of order $5,\ 10$ and 
$16$ in \cite{Berndt3}.

Parnami,Agrawal and Rajwade \cite{Parnami2} showed that for an odd prime $l$, it is sufficient to calculate $J_{l}(1,(l-3)/2)$ number of Jacobi sums 
for $l>3$ and $J_{l}(1,1)$ for $l=3$ to obtain all the Jacobi sums of order $l$. Thus, it reduced the complexity to $l^{2}-(l-3)/2$ for $l>3$ and
$l^{2}-1$ for $l=3$. Acharya and Katre \cite{Acharya1} indicated that calculating all the Jacobi sums of order $2l$ is not essential, and it is 
enough to calculate $J_{2l}(1,n)$ for $1\leq n\leq 2l-3$, $n$ odd or $1\leq n\leq 2l-2$, $n$ even number of Jacobi sums. Ahmed and Tanti \cite{Helal3}
showed that Jacobi sums of order $2l^{2}$ can be determined from the Jacobi sums of order $l^{2}$. The Jacobi sums of order $2l^2$
can also be obtained from $J_{2l^{2}}(1,n)$, $1\leq n\leq 2l^{2}-3$ for $n$ odd (or equivalently, $2\leq n\leq 2l^{2}-2$ for $n$ even). Further the 
Jacobi sums of order $l^{2}$ can be evaluated if one knows the Jacobi sums $J_{l^{2}}(1,i)$, $1\leq i\leq \frac{l^{2}-3}{2}$. 

For some small values of $e$ the study of congruences of Jacobi sums is available in the literature. For $l$ an odd prime
Dickson \cite{Dickson2} obtained the congruences $J_{l}(1,n)\equiv -1 \pmod{(1-\zeta_{l})^{2}}$ for $1\leq n \leq l-1$. Parnami, Agrawal and Rajwade 
\cite{Parnami2} also calculated this separately. Iwasawa \cite{Iwasawa1} in $1975$, and in $1981$ Parnami, Agrawal and Rajwade \cite{Parnami1} 
showed that the above congruences also hold $\pmod{(1-\zeta_{l})^{3}}$. Further in $1995$, Acharya and Katre \cite{Acharya1} extended the work on finding
the congruences for 
Jacobi sums and  showed that 
\begin{center}
$J_{2l}(1,n)\equiv -\zeta_l^{m(n+1)}($mod$\ (1-\zeta_{l})^{2})$,
\end{center} 
where $n$ is an odd integer such that $1\leq n \leq 2l-3$ and $m=$ind$_{\gamma}2$.
Also in $1983$, Katre and Rajwade \cite{Katre2} obtained the congruence of Jacobi sum of order $9$, i.e.,
\begin{center}
$J_{9}(1,1)\equiv -1-($ind$\ 3)(1-\omega)($mod$\ (1-\zeta_{9})^{4})$,
\end{center} 
where $\omega = \zeta_{9}^{3}$.
In $1986$, Ihara \cite{Ihara1} showed that if $k>3$ is an odd prime power, then
 \begin{center}
 $J_{k}(i,j)\equiv -1 \pmod{(1-\zeta_{k})^{3}}$.
 \end{center}
Evans (\cite{Evans1}, $1998$) used simple methods to generalize this result for all $k>2$, getting sharper congruences in some cases, especially 
when $k>8$ is a power of $2$. Congruences for the Jacobi sums of order $l^{2}$ ($l>3$ prime) were obtained by Shirolkar and Katre \cite{Shirolkar1}. 
They showed that\\ 
$J_{l^{2}}(1,n)\equiv
 \begin{cases}
  -1 + \sum_{i=3}^{l}c_{i,n} (\zeta_{l^2} -1)^{i} ($mod$\ (1-\zeta_{l^2})^{l+1}) \ \ \ \ \ if\ $gcd$ (l,n)=1, \\
  -1 \ ($mod$\ (1-\zeta_{l^2})^{l+1}) \ \ \ \ \ \ \  \ \ \ \ \ \ \ \ \ \ \ \ \ \ \ \ \ \ \ \ \ \ \ \ if\ $gcd$ (l,n)=l.  
 \end{cases}$ \\ \\ 
Recently, Ahmed and Tanti \cite{Helal3}(paper is posted in Arxiv and is under the process of publication) determined the congruences $\pmod{(1-\zeta_{l^2})^{l+1}}$ for Jacobi sums of order $2l^{2}$. They split the 
problem into two cases:\\ 
\textbf{Case 1}. $n$ is odd. \\
\textbf{Case 2}. $n$ is even.\\
Further, \textbf{Case 1} splits into three sub-cases:\\
\textbf {Subcase i}. $n=l^{2}$.\\
\textbf {Subcase ii}. $n=dl, \ $where$\ 1\leq d\leq 2l-1, d$ is an odd and $d\neq l$.\\
\textbf {Subcase iii}. gcd$(n,2l^{2})=1$.\\ \\
In the same, they established the relationship of $n$ of Jacobi sums $J_{2l^{2}}(1,n)$, considering $n$ odd as well as even. They showed that, once we 
determined the congruences for $n$ odd, by using the following relation  
\begin{center}
$J_{2l^{2}}(1,n)=\chi_{2l^2}(-1)J_{2l^{2}}(1,2l^{2}-n-1)$,
\end{center}
one can determined the congruences for $n$ even also and conversely. They \cite{Helal3} precisely proved:
\begin{theorem} 
Let $l\geq 3$ be a prime and $q=p^{r}\equiv 1 \pmod {2l^{2}}$. If $1\leq n \leq 2l^{2}-3$ and $1\leq d \leq 2l-1$ are odd integer, then a congruence for $J_{2l^{2}}(1,n)$ over $\mathbb{F}_{q}$ is given by
\begin{equation*}
J_{2l^{2}}(1,n)\equiv \begin{cases}
\zeta_{l^{2}}^{-w}(-1+\sum_{i=3}^{l} c_{i,n}(\zeta_{l^{2}}-1)^{i}) \pmod {(1-\zeta_{l^{2}})^{l+1}},\ \mbox{if $n=l^{2}$}, \\ \\ 
-\zeta_{l^{2}}^{-w(dl+1)}(-1+\sum_{i=3}^{l} c_{i,n}(\zeta_{l^{2}}-1)^{i})(-1+\sum_{i=3}^{l} c_{i,n}(\zeta_{l^{2}}^{(-1-dl)/2}\\-1)^{i})
\pmod {(1-\zeta_{l^{2}})^{l+1}}  \ \mbox{if $d\neq l$ odd integer and $n=dl$},\\ \\
  \zeta_{l^{2}}^{-w(n+1)}(-1+\sum_{i=3}^{l} c_{i,n}(\zeta_{l^{2}}-1)^{i})(-1+\sum_{i=3}^{l} c_{i,n}(\zeta_{l^{2}}^{n}-1)^{i})\\ (-1+\sum_{i=3}^{l} c_{i,n}
  (\zeta_{l^{2}}^{(1-l^{2})/2}-1)^{i})\pmod {(1-\zeta_{l^{2}})^{l+1}},\,\\ \mbox{if \, $gcd(n,2l^{2})=1$},  
\end{cases}
\end{equation*}
where $c_{i,n}$ are as described in the Theorem \ref{T1}.
\end{theorem} 
The congruences of Jacobi sums $J_{9}(1,n)$ of order $9$ which was not fully covered in \cite{Shirolkar1}, and in \cite{Helal3} included the congruences of Jacobi sums $J_{9}(1,n)$ of order $9$ and  revised the result of congruences 
of Jacobi sums of order $l^2$ for $l\geq3$ a prime. Hence Theorem 5.4 \cite{Shirolkar1} was precisely revised as the following 
theorem.
\begin{theorem} \label{T1}
Let $l\geq 3$ be a prime and $p^{r}=q\equiv 1 \pmod {l^{2}}$. If $1\leq n \leq l^{2}-1$, then a congruence for $J_{l^{2}}(1,n)$ for a 
finite field $\mathbb{F}_{q}$ is given by \\ $J_{l^{2}}(1,n)\equiv
\begin{cases}
-1 + \sum_{i=3}^{l}c_{i,n} (\zeta_{l^{2}} -1)^{i} \pmod {(1-\zeta_{l^{2}})^{l+1}} \ \ \ \ \ if\ \emph{gcd} (l,n)=1, \\
 -1 \pmod{(1-\zeta_{l^{2}})^{l+1}} \ \ \ \ \ \ \  \ \ \ \ \ \ \ \ \ \ \ \ \ \ \ \ \ \ \ \ \ \ \ if\ \emph{gcd} (l,n)=l,  
\end{cases}$  \\ 
where for $3\leq i \leq l-1$, $c_{i,n}$ are described by equation (5.3) and $c_{l,n}=S(n)$ is given by Lemma 5.3 in \cite{Shirolkar1}.
\end{theorem}
\section{Cyclotomic number}
Since the time of Gauss, many authors have approached the problem of determining cyclotomic numbers in terms of the solutions of certain Diophantine 
systems. Such a problem arises when Gauss solving his period equation in the case in terms of the uniquely determined $L$ of the Diophantine system 
$4p=L^{2}+27M^{2}$, $L\equiv 1 \pmod 3$. In 1935, a series of three papers were published by Dickson, in which he reviewed and extended the theory of 
cyclotomy. In the first one \cite{Dickson1}, he considered the cases for cyclotomic numbers of order $e\leq 6$, $e=8,\ 10,$ and $12$. In the second paper 
\cite{Dickson2}, he explained the general theory for cyclotomic numbers of prime order and twice of a prime order and he clearly studied the cases 
$e=14$ and $22$. In the third paper \cite{Dickson3}, he discussed cyclotomic numbers of order $e=9$ and $18$, and in the last part of this paper, he 
singled out a part entitles \begin{center}
Theory for $\phi(e)=8,\  e=15,16,20,24,30.$
\end{center} 
However, the case $e = 30$ was completely ignored and one relation associated with $e = 16$ was also omitted. The case $e = 15$ was left with the 
sign of ambiguity and only introductory discussions were given to $e = 20$ and $24$. The sign of ambiguity for $e = 15$ was resolved by Muskat
\cite{Muskat1}. He also provided a complete analysis for $e=24$ and $30$.

The cyclotomic number may be defined for $e=1$; in that case we have $(0,0)_{1}=q-2$. The determination of cyclotomic numbers of order $e=2$ in
$\mathbb{F}_{p}$ was considered by Dickson in \cite{Dickson1} in terms of the period equation [\cite{Dickson1}, \S 9]. He shown that $(a,b)_{2}$'s 
are uniquely determined by $p=ek+1$ and period equation become $\eta^{2}+\eta+c=0$, where $c=kn_{1}-(1,1)_{2}$, if $k$ even $c=-1/4(p-1)$; if $k$ 
odd $c=1/4(p+1)$. Again the case $e=2$ was considered by Whiteman \cite{Whiteman2} in terms of the Jacobsthal sums and his Diophantine system 
become $p=a^{2}+b^{2}$; $a=(\psi_{4}(1)+2)/2$ and $b=\psi_{4}(\gamma)/2$. The delightful book by Davenport \cite{Davenport1}, for $e=2$, the cyclotomic 
numbers do not depend upon the generator $\gamma$, they are given by:
\begin{equation*}
(0,0)_{2}=(q-5)/4,\ \ (0,1)_{2}=(1,0)_{2}=(1,1)_{2}=(q-1)/4, \ \ if \ q\equiv 1 \pmod 4;
\end{equation*}
\begin{equation*}
(0,0)_{2}=(1,0)_{2}=(1,1)_{2}=(q-3)/4,\ \ (0,1)_{2}=(q+1)/4, \ \ if \ q\equiv 3 \pmod 4.
\end{equation*}

The determination of the cyclotomic numbers of order $l=3$ in $\mathbb{F}_{p}$ was considered by Gauss in \cite{Gauss2} in terms of the solutions of 
the diophantine system $4p=L^{2}+27M^{2}$, $L\equiv 1\pmod 3$, when he obtained his period equation in this case in terms of the uniquely determined $L$. 
The three cyclotomic periods of order $3$ satisfy $x^{3}+x^{2}-\dfrac{p-1}{3}x-\frac{1}{27}(3p-1+pl)=0$. These equations determine $L$ uniquely, 
but $M$ is determined only upto sign. Gauss gives formulae for cyclotomic numbers of order $3$ in terms of $L$ and $M$.
\begin{theorem} \cite{Gauss2}
For a prime $p\equiv 1 \pmod 3$, write
\begin{equation*}
4p=L^{2}+27M^{2}, \ L\equiv 1 \pmod 3.
\end{equation*}
Then the $9$ cyclotomic numbers of order $3$ are given by
\begin{equation*}
(0,0)_{3}=(p-8+L)/9,
\end{equation*}
\begin{equation*}
(0,1)_{3}=(1,0)_{3}=(2,2)_{3}=(2p-4-L+9M)/18,
\end{equation*}
\begin{equation*}
(1,1)_{3}=(0,2)_{3}=(2,0)_{3}=(2p-4-L-9M)/18,
\end{equation*}
\begin{equation*}
(1,2)_{3}=(2,1)_{3}=(p+1+L)/9.
\end{equation*}
\end{theorem}
He says that these formulae give the cyclotomic numbers of order $3$ for some generator $\gamma$ of $\mathbb{F}_{p^{*}}$. If $M$ is replaced by $-M$ in 
all the formulae then one gets cyclotomic numbers corresponding to some other generator $\gamma^{\prime}$ of $\mathbb{F}_{p^{*}}$. One says that the 
cyclotomic problem in $\mathbb{F}_{p}$ for $l=3$ was solved by Gauss. However, the solution does not make it clear which sign of $M$ goes with which 
$\gamma$, without an alternative evaluation of some cyclotomic numbers of order $3$, say $(1,0)_{3}$ or $(1,1)_{3}$. In a footnote to the section $358$ 
of \cite{Gauss2} (p. $444$, English edition or p. $432$, German edition) Gauss remarks: \textquotedblleft As far as the ambiguity of the sign of $M$ in 
$4p=L^{2}+27M^{2}$, $L\equiv 1\pmod 3$, for the determination of cyclotomic numbers of order $3$, is concerned, it is unnecessary to consider this 
question here, and by the nature of the case it cannot be determined because it depends on the selection of the primitive root $g$ mod $p$. For some 
primitive roots, $M$ will be positive, for others negative$\textquotedblright$. Later, this case $l=3$ was again taken up by Dickson \cite{Dickson1} 
and considered the same diophantine equation. But his calculation was again a Gauss type of ambiguity. In $1952$, Whiteman \cite{Whiteman2} considered 
the same case and resolved the sign of ambiguity using Jacobsthal sums. He gave the diophantine equation $4p=c^{2}+3d^{2}$; 
$c\equiv 1\pmod 3$ and $d\equiv 0\pmod 3$. Further, M. Hall \cite{Hall1} and Storer \cite{Storer3} generalized the results of Gauss and Dickson for 
$l=3$ to finite fields of $q=p^{r}$ elements. However, when $p\equiv 1\pmod 3$, there results for $\mathbb{F}_{q}$ again have a Gauss and Dickson type 
of ambiguity.

For a prime $p$, the theory for cyclotomy \cite{Dickson1} is well known. The corresponding theory has been developed for $pq$ in \cite{Whiteman1}. 
Further these theories were extended in \cite{Storer2} in connection with the construction of finite difference sets. The above mention cases for $p$ or 
$pq$, the cyclotomic constants depend upon one or more representation as binary quadratic forms of the type $s^{2}+Dt^{2}$. In \cite{Storer3}, T. Storer 
established the uniqueness for $p^{r}$ of cases $e=3,\ 4,\ 6,\ 8$, by same generalizing procedure in \cite{Storer2}. For the cases $e=3$ and $e=4$, 
the unique proper representation are
\begin{equation*}
4p^{r}=s^{2}+27t^{2}; \ \ s\equiv 1\pmod 3
\end{equation*}  
\begin{center}
and
\end{center}
\begin{equation*}
p^{r}=s^{2}+4t^{2}; \ \ s\equiv 1\pmod 4
\end{equation*}
respectively. Cyclotomic numbers of order $e$ over finite $\mathbb{F}_{q}$ with characteristic $p$, in general, can not be determined only in terms of $p$
and $e$, but that one requires a quadratic partition of $q$ too. The cases for $e=6$ and $e=8$, the unique proper representation were established in 
terms of the binary quadratic form $s^{2}+Dt^{2}$ and two such forms respectively.

The formulae for cyclotomic numbers of order $4$ ($p\equiv 1\pmod 4, p$  prime) in terms of the quadratic partition 
$p=s_{0}^{2}+t_{0}^{2},\  s_{0}\equiv 1 \pmod 4$ was obtained by Gauss \cite{Gauss1}, which fixes $t_{0}$ upto sign and $s_{0}$ uniquely. Further, 
Dickson \cite{Dickson1} also worked in this account and his diophantine equation was $p=x^{2}+4y^{2}$, $x\equiv 1\pmod 4$. However, Gauss and Dickson 
did not resolve the sign of ambiguity in $t_{0}$ and $y$ respectively, viz. given a generator $\gamma$ of $\mathbb{F}_{p^{*}}$, it is not clear that
which sign of $t_{0}$ and $y$ gives correct formulae for the cyclotomic numbers corresponding to $\gamma$. Again the case $e=4$ was considered by 
Whiteman \cite{Whiteman2} and resolve the sign of ambiguity using Jacobsthal sums. The corresponding result of M. Hall \cite{Hall1} for $\mathbb{F}_{q}$ 
by setup $q=p^{r}\equiv 1 \pmod 4$ also has a similar sign of ambiguity in the case when $q=p\equiv 1 \pmod 4$. Later Katre and Rajwade \cite{Katre1} 
resolve the sign of ambiguity and they gave the formulae to determined the cyclotomic numbers of order $4$ as:\\
for $k$ even, 
\begin{align*}
(0,0)_{4}&=1/16(q-11-6s),\\
(0,1)_{4}&=1/16(q-3+2s+4t),\\
(0,2)_{4}&=1/16(q-3+2s),\\
(0,3)_{4}&=1/16(q-3+2s-4t),\\
(1,2)_{4}&=1/16(q+1-2s),\\
\end{align*}    
and for $k$ odd,
\begin{align*}
(0,0)_{4}&=1/16(q-7+2s),\\
(0,1)_{4}&=1/16(q+1+2s-4t),\\
(0,2)_{4}&=1/16(q+1-6s),\\
(0,3)_{4}&=1/16(q+1+2s+4t),\\
(1,1)_{4}&=1/16(q-3-2s).\\
\end{align*}

The next case $l=5$ was treated by Dickson \cite{Dickson1} for prime $p$, using the properties of Jacobi sums. Dickson considers the Diophantine equation
$16p=x^{2}+50u^{2}+50v^{2}+125w^{2}$, $v^{2}-4uv-u^{2}=xw$, and $x\equiv 1 \pmod 5$. These system has exactly eight integral simultaneous solutions. 
If $(x,u,v,w)$ is one solution, also $(x,-u,-v,w)$ and $(x,\pm u,\mp v,-w)$ are solutions. The remaining four are derived from these four by changing all
signs. In terms of these solutions Dickson gives formulae for cyclotomic numbers of order $5$. Again, Dickson does not tell which solution goes with 
which $\gamma$. Later, Whiteman \cite{Whiteman2} considered the same diophantine system for cyclotomic numbers of order $5$ and resolve the sign of 
ambiguity using Jacobsthal sums. Katre and Rajwade \cite{Katre4} for the determination of a unique solution, they considered a fifth root of unity in 
terms of a solution of Dickson's diophantine system.

Determination of cyclotomic numbers of order $e=6$, Dickson \cite{Dickson1} showed that $36$ cyclotomic constants $(a,b)_{6}$ depend solely upon the 
decomposition $A^{2}+3B^{2}$ of the prime $p=6k+1$. Due to the signs of ambiguity, he took $2$ be a cubic residue of $p$, $\gamma^{m}\equiv 2 \pmod p$ 
and $m\equiv 1\ or \ 4\pmod 6$.

The case $e=7$ \cite{Leonard1}, the cyclotomic numbers can be given in terms of Dickson-Hurwitz sums using the work of Muskat \cite{Muskat2} or a theorem 
of Whiteman \cite{Whiteman3}. In \cite{Leonard1}, Leonard and Williams obtained the cyclotomic numbers of order $e=7$ in terms of the solutions of a 
certain triple of Diophantine equations, analogous to the expressions for the cyclotomic numbers of order $5$ in terms of the solutions of a pair of 
Diophantine equations \cite{Whiteman3}. They used the result of Muskat \cite{Muskat2} to evaluate the cyclotomic numbers of order $7$. If 
$p\equiv 1\pmod 7$ then there are exactly six integral simultaneous solutions of the triple of Diophantine equations
\begin{equation*}
2x_{1}^{2}+42(x_{2}^{2}+x_{3}^{2}+x_{4}^{2})+343(x_{5}^{2}+3x_{6}^{2})=72p,
\end{equation*}
\begin{equation*}
12x_{2}^{2}-12x_{4}^{2}+147x_{5}^{2}-441x_{6}^{2}+56x_{1}x_{6}+24x_{2}x_{3}-24x_{2}x_{4}+48x_{3}x_{4}+98x_{5}x_{6}=0,
\end{equation*}
\begin{align*}
& 12x_{3}^{2}-12x_{4}^{2}+49x_{5}^{2}-147x_{6}^{2}+28x_{1}x_{5}+28x_{1}x_{6}+48x_{2}x_{3}+24x_{2}x_{4}+24x_{3}x_{4}\\ &+490x_{5}x_{6}=0,
\end{align*}
satisfying $x_{1}\equiv 1\pmod 7$, distinct from the two trivial solutions $(-6t,$ $\pm 2u,$ $\pm 2u,$ $\mp 2u,$ $0,$ $0)$, where $t$ is given uniquely 
and $u$ is given ambiguously by
\begin{equation*}
p=t^{2}+7u^{2}, \ \ \ \ \ \ t\equiv 1 \pmod 7.
\end{equation*}
If $(x_{1},x_{2},x_{3},x_{4},x_{5},x_{6})$ is a nontrivial solution with $x_{1}\equiv 1\pmod 7$ then two others solutions are given by 
$(x_{1},-x_{3},x_{4},x_{2},(-x_{5}-3x_{6})/2,(x_{5}-x_{6})/2)$ and $(x_{1},-x_{4},x_{2},-x_{3},(-x_{5}+3x_{6})/2,(-x_{5}-x_{6})/2)$. Each of the other 
three can be obtained from one given above by changing the signs of $x_{2},x_{3},x_{4}$. The result obtained in \cite{Leonard1} is almost similar for 
$p\equiv 1\pmod 5$, result obtained in \cite{Dickson1}, and which is implicit in the work of Dickson \cite{Dickson1},\cite{Dickson2}, does not appear 
in the literature.

Dickson \cite{Dickson1} showed in the case $e=8$ that $64$ cyclotomic constants $(a,b)_{8}$ depend solely upon the decompositions 
$p=x^{2}+4y^{2}$ and $p=a^{2}+2b^{2}$; $x\equiv a \ \equiv 1\pmod 4$, where the signs of $y$ and $b$ depend on the choice of the generator $\gamma$. 
There are four sets of formulaes depending on whether $k$ is even or odd and whether $2$ is a biquadratic residue or not. Further, Lehmer \cite{Lehmer1} 
improve the results of Dickson and gave the complete table of cyclotomic numbers $(a,b)_{8}$ of order $8$.

The cyclotomic problem for $e=9$ was studied by Dickson \cite{Dickson3} and he gave a simple complete theory for $e=9$. Each cyclotomic numbers are 
expressed as a constant plus a linear combination of $p,\ L,\ M,\ c_{0},\ c_{1},\ c_{2},\ c_{3},\ c_{4},\ c_{5}$ where 
\begin{equation*}
4p=L^{2}+27M^{2}, \ \ \ \ \ L\equiv 7\pmod 9
\end{equation*} 
and 
\begin{center}
$p=(\sum_{i=0}^{5}c_{i}\beta^{i})(\sum_{i=0}^{5}c_{i}\beta^{-i})$\ \ \ \ ($\beta$ \ be \ a\ primitive\ ninth\ root\ of \ unity)
\end{center}
is a factorization of $p$ in the field of ninth roots of unity. Further, he gave theorem to choose $M$. In $1967$, again the case $e=9$ were considered 
by Baumert and Fredricksen \cite{Baumert1}. They carried out the result of Dickson \cite{Dickson3}. But they worked little bit further
Dickson \cite{Dickson3} did. They gave simple relation to choose $M$ as Dickson gave. 

The case $e=11$ were considered by Leonard and Williams \cite{Leonard2}. They considered the Diophantine equations
\begin{equation*}
1200p=12w_{1}^{2}+33w_{2}^{2}+55w_{3}^{2}+110w_{4}^{2}+330w_{5}^{2}+660(w_{6}^{2}+w_{7}^{2}+w_{8}^{2}+w_{9}^{2}+w_{10}^{2}),
\end{equation*} 
\begin{align*}
0&=45w_{2}^{2}+5w_{3}^{2}+20w_{4}^{2}-540w_{5}^{2}+720w_{6}^{2}-720w_{10}^{2}-288w_{1}w_{5}+30w_{2}w_{3}\\ &\ \ \ \-120w_{2}w_{4}-72w_{2}w_{5}  +200w_{3}w_{4}-360w_{3}w_{5}+360w_{4}w_{5}+1440w_{6}w_{7}\\ &\ \ \ \-1440w_{6}w_{8}+1440w_{7}w_{8}-1440w_{7}w_{9} +1440w_{8}w_{9}-1440w_{8}w_{10}+2880w_{9}w_{10},
\end{align*}
\begin{align*}
0&=45w_{2}^{2}-35w_{3}^{2}-80w_{4}^{2}+720w_{9}^{2}-720w_{10}^{2}-144w_{1}w_{4}-144w_{1}w_{5} \\ &\ \ \ \ +150w_{2}w_{3}-96w_{2}w_{4} -216w_{2}w_{5}+160w_{3}w_{4}+120w_{3}w_{5}+240w_{4}w_{5}\\ &\ \ \ \ +2880w_{6}w_{7}-1440w_{6}w_{9}+1440w_{7}w_{8} -1440w_{7}w_{10}+1440w_{8}w_{9}\\ &\ \ \ \ +1440w_{8}w_{10}+1440w_{9}w_{10},
\end{align*}
\begin{align*}
0&=45w_{2}^{2}+5w_{3}^{2}+20w_{4}^{2}-540w_{5}^{2}+720w_{7}^{2}-720w_{10}^{2}-96w_{1}w_{3}-48w_{1}w_{4}\\ &\ \ \ \-144w_{1}w_{5}+126w_{2}w_{3}  +108w_{2}w_{4}-360w_{2}w_{5}+20w_{3}w_{4}-60w_{3}w_{5}\\ &\ \ \ \ +600w_{4}w_{5}+1440w_{6}w_{7}+1440w_{6}w_{8} -1440w_{6}w_{10} +1440w_{7}w_{8}\\ &\ \ \ \ +1440w_{7}w_{10}+1440w_{9}w_{10}+2880w_{8}w_{9},
\end{align*}
\begin{align*}
0&=27w_{2}^{2}+35w_{3}^{2}-40w_{4}^{2}-360w_{5}^{2}+720w_{8}^{2}-720w_{10}^{2}-72w_{1}w_{2}-24w_{1}w_{3}\\ &\ \ \ \ -48w_{1}w_{4}-144w_{1}w_{5}  +114w_{2}w_{3}+48w_{2}w_{4}+144w_{2}w_{5}+320w_{3}w_{4} \\ &\ \ \ \ +1440w_{6}w_{7}+1440w_{6}w_{9}+1440w_{6}w_{10} +2880w_{7}w_{8}+1440w_{7}w_{9}\\ &\ \ \ \ +1440w_{8}w_{9}+1440w_{9}w_{10},
\end{align*}
\begin{equation*}
w_{3}+2w_{4}+2w_{5}\equiv 0\pmod {11},
\end{equation*}
\begin{equation*}
w_{2}-w_{4}+3w_{5}\equiv 0\pmod {11}.
\end{equation*}
In terms of the solutions of above Diophantine euations, they gave the complete formulae for cyclotomic numbers of order $11$. Further, determination of 
cyclotomic problem is somewhat incomplete, when they finding the solutions of above mentioned Diophantine systems using Jacibsthal-Whiteman sums. 

The cycloyomic problems for $e=12$ was considered by Dickson \cite{Dickson1} and he showed that cyclotomic constants $(a,b)_{12}$ depend solely upon the 
decomposition $p=x^{2}+4y^{2}$ and $p=A^{2}+3B^{2}$ of the prime $p=12k+1$, where $x\equiv 1 \pmod 4$ and $A\equiv 1 \pmod 6$. But his analysis depends 
upon elaborate computations and is not entirely definitive. For settlement of signs of ambiguity for case $e=12$, he considered lots of small cases, viz. 
$2$ be a cubic residue of $p$, $3$ be a biquadratic residue and non-residue of $p$. For odd case, again he considered $2m\equiv 2 \pmod {12}$ and 
$2m\equiv 10 \pmod {12}$. For even case $2m\equiv 8 \pmod {12}$ and $2m\equiv 4 \pmod {12}$. Later, the same case was considered by Whiteman
\cite{Whiteman4} in a different direction. To evaluate the complete solution of cyclotomic constants, he divided the prime into $12$ different classes 
and obtained formulae holding for different classes.

The study of cyclotomic numbers of order $14$ started by Dickson \cite{Dickson2} in $1935$. Dickson proved that it is possible to represent the cyclotomic
numbers $(a,b)_{14}$ as a linear combination of the $B_{e}(i,v)$, if $e$ is an odd prime or twice of an odd prime. Dickson's result that, given the septic
character of $2$, the $B_{14}(i,v)$ can be given as a linear combination of the $B_{7}(i,v)$ is then employed to express the $(a,b)_{14}$ in terms of the 
$B_{7}(i,v)$. Muskat \cite{Muskat2} determined the cyclotomic numbers $(a,b)_{14}$ explicitly in terms of the $B_{e}(i,v)$ where $e$ is twice of an odd 
prime and the transformation is due to the Whiteman \cite{Whiteman3}. 

The study of cyclotomic numbers of order $15$ was began by Dickson \cite{Dickson3} in $1935$ and completed by Muskat \cite{Muskat1} in $1968$. But Dickson
calculation have a sign of ambiguity. Muskat \cite{Muskat1} resolve the sign of ambiguity. In $1986$, Frisen, Muskat, Spearman and Williams \cite{Friesen1}
considered the same case by setup $q=p^{2}\equiv 1\pmod {15}$. They showed that as $p\equiv 4\pmod {15}$, $28$ different cyclotomic numbers were evaluated
by
\begin{align*}
p&=A^{2}-AB+B^{2},\ \ A\equiv-1 \pmod 3,\ \ B\equiv 0 \pmod 3,\\
p& = T^{2}+15U^{2},\ \ T\equiv-1 \pmod 3,
\end{align*} 
and $p\equiv 11\pmod {15}$, $29$ different cyclotomic numbers were evaluated by
\begin{align*}
p&=X^{2}+5U^{2}+5V^{2}+5W^{2},\ \ X\equiv-1 \pmod 5,\\
XW&= V^{2}-UV-U^{2}.
\end{align*}
In \cite{Buck1}, Buck, Smith, Spearman and Williams used Dickson and Muskat evaluations of the Jacobi sums of order $15$ to obtain the values of the 
Dickson-Hurwitz $B_{15}(i,v)$ of order $15$ defined by
\begin{equation*}
B_{15}(i,v)=\sum_{h=0}^{14}(h,i-vh)_{15}.
\end{equation*}
Further they use a special case of a theorem of Friesen, Muskat, Spearman and Williams \cite{Friesen1} is used to express each cyclotomic number in 
terms of the Dickson-Hurwitz sums and then using the values for the Dickson-Hurwitz sum, they derived an explicit formulae for the cyclotomic numbers of 
order $15$. Each cyclotomic numbers of order $15$ can be expressed as an integral linear combination of the integers 
$p,$ $1,$ $a,$ $b,$ $c,$ $d,$ $x,$ $u,$ $v,$ $w,$ $b_{0}$ $b_{1}$ $b_{2}$ $b_{3}$ $b_{4}$ $b_{5}$ $b_{6}$ $b_{7}$. The integers 
$a,$ $b,$ $ c,$ $d,$ $x,$ $u,$ $v,$ $w$ have the following properties:
\begin{align*}
p=a^{2}+3b^{2}, \ \ \ \ a\equiv -1\pmod 3,\\ 
p=c^{2}+15d^{2}, \ \ \ \ c\equiv -1\pmod 3,\\ 
p=x^{2}+5u^{2}+5v^{2}+5w^{2},\\
xw=v^{2}-uv-u^{2}, \ \ \ \ x\equiv -1\pmod 5.
\end{align*}

E. Lehmar \cite{Lehmer2} raised the question whether or not constants $\alpha,\ \beta, \ \gamma,\ \delta, \ \epsilon$ can be found such that
\begin{equation}\label{lehmer1.1}
265(a,b)_{16}=p+\alpha x+\beta y+\gamma a+\delta b+\epsilon,
\end{equation}
at least for some $(a,b)_{16}$ after written the article \cite{Lehmer1}. To answer this question, she undertook
 the following experiment on the SWAC (National Bureau of Standards Western Automatic Computer). The cyclotomic constants of order sixteen were computed 
 for eight primes $p$ of the form $32n+1$ for which $2$ is not a biquadratic residue. She found that (\ref{lehmer1.1}) is not satisfied for any 
 $(a,b)_{16}$ when the signs of $y$ and $b$ are taken in accordance with the results on cyclotomic constants of order eight. The calculations exhibited 
 eight solutions, while the formula gave only six. A similar computation for primes $p$ of the form $32n+17$ also led to a negative result. She 
 regretfully conclude that the cyclotomic constants of order sixteen are not expressible in terms of these \cite{Dickson3} quadratic partitions alone. 
 The SWAC experiment leaves open the question of determining if the equation (\ref{lehmer1.1}) can be satisfied for any prime $p$ for which $2$ is a 
 biquadratic residue. In \cite{Whiteman5}, Whitwman gave formulae for cyclotomic constants of order sixteen in affirmative for six of the cyclotomic 
 constants in terms of parameters of quartic, octic and bioctic Jacobi sums. He further gave a table of formulas for $(a,0)_{16}$. The following formulas
 involving only $p,\ a$ and $x$ are derived. Let $p=16k+1$ be a prime. If the integer $m$ is defined by the congruence $\gamma^{m}\equiv 2\pmod p$, then
 \begin{align*}
 256(0,0)_{16}&=p-47-18x\ \ \ \ \ \ \ \ \ \ \ (k \ even, m\equiv 4\pmod 8),\\
 256(8,0)_{16}&=p-15-18x-32a\ \ \ (k \ even, m\equiv 4\pmod 8),\\
 256(4,8)_{16}&=p+1-2x\ \ \ \ \ \ \ \ \ \ \ \ \ \ (k \ even, m\equiv 0\pmod 8),\\
 256(0,0)_{16}&=p-31-18x-16a\ \ \ (k \ odd, m\equiv 0\ \ \pmod 8),\\
 256(0,8)_{16}&=p+1-18x-48a\ \ \ \ \ (k \ odd, m\equiv 0\ \ \pmod 8),\\
 256(4,0)_{16}&=p-15-2x+16a\ \ \ \ \ (k \ odd, m\equiv 4\ \ \pmod 8),
 \end{align*}
 where the signs of $a$ and $x$ are selected so that $a\equiv x\equiv 1\pmod 4$. In other cases it is shown that the cyclotomic constants $(a,b)_{16}$
 are such that $256(a,b)_{16}$ is expressible as a linear combination with integer coefficients of $p,$ $a,$ $b,$ $x,$ $y$ and certain other integers
 $c_{0},$ $c_{1},$ $c_{2},$ $c_{3},$ $d_{0},$ $d_{1},$ $d_{2},$ $d_{3},$ $d_{4},$ $d_{5},$ $d_{6},$ $d_{7}$. Evans and Hill \cite{Evans3} gave complete
 table for cyclotomic numbers of order sixteen. The computations were performed on the Burroughs $6700$ at UCSD, with use of algorithms described 
 in \cite{Whiteman5}. It is not possible to accomplish sign resolutions with the use of formulas from Whiteman's, so it has been utilized \cite{Evans2}
 to gave elementary resolutions of sign ambiguities in quartic and octic Jacobi and Jacobsthal sums, in certain cases. This is what motivated the worked 
 of Evans and Hill \cite{Evans3}.
 
Dickson \cite{Dickson3} gave relation to determine the cyclotomic numbers of order $18$. But he did not provide a complete table for cyclotomic numbers
of order $18$. Baumert and Fredricksen \cite{Baumert1} considered the case of cyclotomic numbers of order $18$ again. They splitting the solution into
cases and introducing the parameters $B=Ind\ 2$ and $T=Ind\ 3$. Introduced parameters reduced actual parameters to 
$p,\ L,\ M,\ c_{0},\ c_{1},\ c_{2},\ c_{3},\ c_{4},\ c_{5}$ as appeared in the determination of the cyclotomic numbers of order $9$. 
They gave complete listing of all  cyclotomic numbers of order $18$ ($k$ odd and $k$ even) in an unpublished mathematical tables file of Mathematics 
of Computation. Further, by the use of cyclotomic formulas given in Table $5$ and $6$, they proved theorem \ref{thm5.2} as an application to difference 
sets.
\begin{theorem} \label{thm5.2}
The only residue difference set or modified residue difference set which exists for $e=18$ is the trivial $19-1-0$ difference set.
\end{theorem}

In \cite{Dickson3}, Dickson gave a sketchy discussions for cyclotomic numbers of order $e=20$. He did not gave the exact formulae for cyclotomic numbers 
of order twenty. In $1970$, Muskat and Whiteman \cite{Muskat3} gave the complete formulae for cyclotomic numbers of order $e=20$. They obtained the 
cyclotomic numbers of order $20$ by the values of the appropriate cyclotomic numbers of order two and four and Jacobi sums of order $5,\ 10$ and $20$.

In \cite{Dickson2}, Dickson completely analyzed the Jacobi sums of order double of a prime. By using those relation, he gave complete determination of 
cyclotomic numbers of order double of an odd prime. Further, in the same paper he gave the general theory for diophantine equations when $e$ be an odd 
prime.

In $1982$, PAR \cite{Parnami2} considered the general $e$ cases; $e$ be an odd prime, $q=p^{r}$, $p\equiv 1\pmod e$. They indicate a general method for 
solving the cyclotomic problem over $\mathbb{F}_{q}$. They calculated the cyclotomic numbers upto $e\leq 19$, onwards the class number of 
$\mathbb{Q}(2 \pi i/e)>1$. They considered the Diophantine systems (slightly different formulation)
\begin{equation*}
q=\sum_{i=0}^{e-1}a_{i}^{2}\sum_{i=0}^{e-1}a_{i}a_{i+1},\ (i.e.,\ 2q=(a_{0}-a_{1})^{2}+(a_{1}-a_{2})^{2}+...+(a_{e-1}-a_{0})^{2}),
\end{equation*}
\begin{equation*}
\sum_{i=0}^{e-1}a_{i}a_{i+1}=\sum_{i=0}^{e-1}a_{i}a_{i+2}=...=\sum_{i=0}^{e-1}a_{i}a_{i+(e-1)/2},
\end{equation*}
\begin{equation*}
1+a_{0}+a_{1}+...+a_{e-1}\equiv 0\pmod e,
\end{equation*}
\begin{equation*}
a_{1}+2a_{2}+3a_{3}+...+(e-1)a_{e-1}\equiv 0\pmod e,
\end{equation*}
which generalizes the diophantine systems of Gauss-Dickson and Leonard-Williams. Moreover they gave a rejection condition 
\begin{equation*}
p\nmid \prod_{\lambda ((n+1)a)>a} H^{\sigma_{a}},
\end{equation*}
which fixes certain Jacobi sums upto conjugate. They gave full details of the transformation (for $r=1$) connecting above mentioned diophantine equation 
with the classical ones for $l=3,\ 5$, but they did not connect the cases for $e>5$ because the rejection condition was too complicated. They again did 
not resolve the signs of ambiguity of Gauss and Dickson type.

Further in $1985$, Katre and Rajwade \cite{Katre3} solved the cyclotomic problem for any prime $e$ in $\mathbb{F}_{q}$, $q=p^{r}$, $p\equiv 1\pmod e$. 
For solving the cyclotomic problem, they added a new condition to Parnami, Agrawal and Rajwade \cite{Parnami2} diophantine systems, i.e. 
$p|\overline{H}\prod_{\lambda((n+1)a)>a}(b-\zeta_{l}^{\sigma_{a^{-1}}})$, where $a^{-1}$ was taken mod $l$, $l$ be an odd prime and $(a,l)=1$. 
For cyclotomic numbers of order $2$ in $\mathbb{F}_{q}$, they told that one can determined the cyclotomic numbers by 
$(0,0)_{2}=(q-5)/4$, $(0,1)_{2}=(1,0)_{2}=(1,1)_{2}=(q-1)/4$ if $q\equiv 1\pmod 4$, $(0,0)_{2}=(1,0)_{2}=(1,1)_{2}=(q-3)/4$, $(0,1)_{2}=(q+1)/4$ if
$q\equiv 3\pmod 4$. The problem arises in Parnami, Agrawal and Rajwade \cite{Parnami2} to connect the Jacobi sums and $e^{th}$ root of unity, 
they considered $\gamma^{(q-1)/e}\equiv 1\pmod e$ gave the proper connection between the Jacobi sums and the $e^{th}$ root of unity mod $p$. 
Additional condition in Katre and Rajwade \cite{Katre3} also resolve the sign of ambiguity arise in Parnami, Agrawal and Rajwade \cite{Parnami2}.

For $l$ an odd prime, Acharya and Katre \cite{Acharya1} determined the cyclotomic numbers of order $2l$ over the field $\mathbb{F}_{q}$ for $q=p^{r}$ 
with the prime $p\equiv 1 \pmod {2l}$ in terms of the solutions of the diophantine systems considered for the $l$ case except that the proper choice of
the solutions for the $2l$ case was made by additional conditions $(iv^{\prime})$, $(v^{\prime})$, $(vi^{\prime})$ which replace the conditions $(iv)$,
$(v)$, $(vi)$ used in the $l$ case. These additional conditions determine required unique solutions thereby also giving arithmetic characterization of 
the relevant Jacobi sums and then the cyclotomic numbers of order $2l$ are determined unambiguously by the following theorem. In the same, they also 
showed how the cyclotomic numbers of order $l$ and $2l$ can be treated simultaneously.
\begin{theorem}\cite{Acharya1}
Let $p$ and $l$ be odd rational primes, $p\equiv 1 \pmod l$ (thus $p\equiv 1 \pmod {2l}$ also), $q=p^{r}$, $r\geq 1$. Let $q=2lk+1$. Let $\zeta_{l}$ 
and $\zeta_{2l}$ be fixed primitive $l-th$ and $2l-th$ roots of unity respectively. Let $\zeta_{l}$ and $\zeta_{2l}$ be related by 
$\zeta_{l}=\zeta_{2l}^{2}$, i.e. $\zeta_{2l}=-\zeta_{l}^{(l+1)/2}$. Let $\gamma$ be a generator of $\mathbb{F}_{q}^{*}$. Let $b$ be a rational integer 
such that $b=\gamma^{(q-1)/l}$ in $\mathbb{F}_{p}$. Let $m=ind_{\gamma}2$. Let $J_{l}(i,j)$ and $J_{2l}(i,j)$ denote the Jacobi sums in 
$\mathbb{F}_{q}$ of order $l$ and $2l$ (respectively)related to $\zeta_{l}$ and $\zeta_{2l}$ (respectively). For $(m,l)=1$, let $\sigma_{m}$ denote 
the automorphism $\zeta_{l}\mapsto \zeta_{l}^{m}$ of $\mathbb{Q}(\zeta_{l})$. For $(m,2l)=1$, let $\tau_{m}$ denote the automorphism 
$\zeta_{2l}\mapsto \zeta_{2l}^{m}$ of $\mathbb{Q}(\zeta_{l})$. Thus if $m$ is odd then $\sigma_{m}=\tau_{m}$ and if $m$ is even then 
$\sigma_{m}=\tau_{m+l}$. Let $\lambda(r)$ and $\varLambda(r)$ denote the least non-negative residues of $r$ modulo $l$ and $2l$ (resp.). 
Let $a_{0}$, $a_{1}$, $\dots$, $a_{l-1}\in \mathbb{Z}$ and let $H=\sum_{i=0}^{l-1}a_{i}\zeta_{l}^{i}$. Consider the arithmetic conditions 
(or diophantine system)\\
$(i)\ \ q=\sum_{i=0}^{l-1}a_{i}^{2}\sum_{i=0}^{l-1}a_{i}a_{i+1},$\\
$(ii)\ \ \sum_{i=0}^{l-1}a_{i}a_{i+1}=\sum_{i=0}^{l-1}a_{i}a_{i+2}=...=\sum_{i=0}^{l-1}a_{i}a_{i+(e-1)/2},
$\\ $
(iii)\ \ 1+a_{0}+a_{1}+...+a_{e-1}\equiv 0\pmod l.$ \\
Let $1\leq n \leq l-2$. If $a_{0}$, $a_{1}$, $\dots$, $a_{l-1}$ satisfy $(i)-(iii)$ together with the additional conditions\\
$(iv)\ \ a_{1}+2a_{2}+3a_{3}+...+(l-1)a_{l-1}\equiv 0\pmod l,$
\\
$(v)\ \ p\nmid \prod_{\lambda ((n+1)m)>m} H^{\sigma_{m}},$
\\
$(vi)\ \ p|\overline{H}\prod_{\lambda((n+1)m)>m}(b-\zeta_{l}^{\sigma_{m^{-1}}})$, where $m^{-1}$ was taken $\pmod l$, \\
then $H=J_{l}(1,n)$ for this $\gamma$ and conversely.
\\
Let $1\leq n \leq 2l-3$ be an odd integer. If $a_{0}$, $a_{1}$, $\dots$, $a_{l-1}$ satisfy $(i)-(iii)$ together with the additional conditions\\
$(iv)^{\prime}\ \ a_{1}+2a_{2}+3a_{3}+...+(l-1)a_{l-1}\equiv u(n+1)\pmod l,$
\\
$(v)^{\prime}\ \ p\nmid \prod_{\varLambda ((n+1)m)>m} H^{\tau_{m}},$
\\
$(vi)^{\prime}\ \ p|\overline{H}\prod_{\varLambda((n+1)m)>m}(b-\zeta_{l}^{\tau_{m^{-1}}})$, where $m^{-1}$ was taken $\pmod {2l}$, \\
then $H=J_{2l}(1,n)$ for this $\gamma$ and conversely.\\
(In $(v)^{\prime}$ and $(vi)^{\prime}$, $m$ varies over only those values which satisfy $1\leq m \leq 2l-1$ and $(m,2l)=1$.)

Moreover, for $1\leq n \leq l-2$ if $a_{0}$, $a_{1}$, $\dots$, $a_{l-1}$ satisfy the conditions $(i)-(vi)$ and if we fix $a_{0}=0$ at the outset and 
write the $a_{i}$ corresponding to a given $n$ as $a_{i}(n)$ then we have $J_{l}(1,n)=\sum_{i=1}^{l-1}a_{i}(n)\zeta_{l}^{i}$ and the cyclotomic numbers 
of order $l$ are given by:\\
$l^{2}(i,j)_{l}=q-3l+1+\varepsilon(i)+\varepsilon(j)+\varepsilon(i-j)+l\sum_{n=1}^{l-2}a_{in+j}(n)-\sum_{n=1}^{l-2}\sum_{k=1}^{l-1}a_{k}(n)$ \\
where
\begin{center}
$\varepsilon(i)= \begin{cases}
0 \hspace*{8mm} \text{if } l|i,\\
l \hspace*{8mm} otherwise,
\end{cases}$
\end{center}
and the subscripts in $a_{in+j}(n)$ are considered modulo $l$.

Similarly, for $n$ odd, $1\leq n \leq 2l-3$, if $a_{0}$, $a_{1}$, $\dots$, $a_{l-1}$ satisfy the conditions $(i)-(iii)$ and $(iv)^{\prime}-(vi)^{\prime}$ 

and if we fix $a_{0}=0$ at the outset and write the $a_{i}$ corresponding to a given $n$ as $b_{i}(n)$ then we have 
$J_{2l}(1,n)=\sum_{i=1}^{l-1}b_{i}(n)\zeta_{l}^{i}$ and the $4l^{2}$ cyclotomic numbers $(i,j)_{2l}$ are given by:\\ 
\begin{align*}
&4l^{2}(i,j)_{2l} = q-3l+1+\varepsilon(i)+\varepsilon(j)+\varepsilon(i-j)+l\sum_{n=1}^{l-2}a_{in+j}(n)-\sum_{n=1}^{l-2}\sum_{k=1}^{l-1}a_{k}(n) \\ & -\{(-1)^{j}+(-1)^{i+k}+(-1)^{i+j}\} \{l+\sum_{k=0}^{l-1}b_{k}(l) +\sum_{u=0}^{l-2}\sum_{k=0}^{l-1}b_{k}(2u+1)\}\\
&+(-1)^{j}l\{b_{\nu(-i)}(l) +\sum_{u=0}^{l-1}b_{\nu(j-2iu-2i)}(2u+1)\}+(-1)^{i+j}l\{b_{\nu(j)}(l) \\
&  +\sum_{u=0}^{l-1}b_{\nu(i+2ju+j)}(2u+1)\}+(-1)^{i+k}l\{b_{\nu(-j)}(l) +\sum_{u=0}^{l-1}b_{\nu(i-2ju-2j)}(2u+1)\},
\end{align*}
where
\begin{center}
$\nu(j)= \begin{cases}
\varLambda(j)/2 \hspace*{8mm} \text{if } j \text{ is even},\\
\varLambda(j+l)/2 \hspace*{2mm} \text{if } j \text{ is  odd},
\end{cases}$
\end{center}
and $\varLambda(r)$ is defined as the least non-negative residue of $r$ modulo $2l$.
\end{theorem}
For $q \equiv 1 \pmod l$, but $p$ does not necessary $\equiv 1 \pmod l$. Let $r$ be the least positive integer such that $q=p^{r}\equiv 1 \pmod l$. In 
the case when $r$ is even, Anuradha and Katre \cite{Anuradha1} obtained the Jacobi sums and cyclotomic numbers of order $l$ and $2l$ just in terms of $q$. 
For $e=l,\ 2l$ and when $r$ is odd, the cyclotomic problem was treated by Anuradha \cite{Anuradha2} in her Ph.D. thesis, and thus the problem for the order
$e$ is settled for all $q \equiv 1 \pmod e$, for $e=l,\ 2l$. Again, for such primes, Shirolkar and Katre \cite{Shirolkar1} obtained formulae for cyclotomic
numbers of order $l^{2}$ in terms of the coefficients of Jacobi sums of orders $l$ and $l^{2}$.
In the same line, recently, Ahmed, Tanti and Hoque \cite{Helal1} (the paper is posted in Arxiv and is under the process of publication) determine an 
expression for determination of cyclotomic numbers of order $2l^{2}$ in 
terms of the coefficient of Jacobi sums of order $l$, $2l$, $l^{2}$ and $2l^{2}$.
\section{Concluding remarks}
Recently, Ahmed, Tanti, Chakraborty and Pusph \cite{Helal2} (The paper is under the process of publication) gave fast computational algorithm to obtain 
all the cyclotomic numbers of order $2l^{2}$. 
They also introduce a new idea to construct matrix of cyclotomic numbers i.e. cyclotomic matrix. Initially, in \cite{Helal2} we have showed that for 
$l=3$, it is enough to calculate $64$ distinct cyclotomic numbers of order $2l^2$ and for $l\neq 3$, it is sufficient to calculate $2l^2+{(2l^2-1)(2l^2-2)}/6$ 
distinct cyclotomic numbers of order $2l^2$. For showing this, they developed an algorithm for equality of cyclotomic numbers. Further, they employ 
those matrices of order $2l^{2}$ and developed a public key cryptosystem. Earlier, K. H. Leung, S. L. Ma and B. Schmidt \cite{Leung1} constructed 
Hadamard matrices of order $4p^{2}$ obtained from Jacobi sums of order $16$ and they proved the following result.
\begin{theorem} \cite{Leung1}
Let $p \equiv 7 \pmod {16}$ be a prime. Then there are integers $a$, $b$, $c$, $d$ with
\begin{equation*}
a \equiv 15 \pmod {16},
\end{equation*}
\begin{equation*}
b \equiv 0 \pmod 4,
\end{equation*}
\begin{equation*}
q^{2}= a^{2}+2(b^{2}+c^{2}+d^{2}),
\end{equation*}
\begin{equation*}
2ab = c^{2}-2cd-d^{2}.
\end{equation*}
If
\begin{equation*}
q= a \pm 2b \ \ or
\end{equation*}
\begin{equation*}
q= a+\delta_{1}b+4\delta_{2}c+4\delta_{1}\delta_{2}4d \ \ with \ \ \delta_{i}=\pm 1,
\end{equation*}
then there is a regular Hadamard matrix of order $4q^{2}$.
\end{theorem}
Betsumiya, Hirasaka, Komatsu and Munemasa \cite{Betsumiya1} gave upper bounds for cyclotomic numbers of order $e$ over a finite field with $q$ elements, where $e$ is a positive divisor of $q-1$. In particular, they showed that under certain assumptions, cyclotomic numbers are at most $\Big \lceil \frac{k}{2}\Big \rceil$, and the cyclotomic number $(0,0)_{e}$ is at most $\Big \lceil \frac{k}{2}\Big \rceil-1$, where $k=(q-1)/e$. They proved the following:
\begin{theorem}\cite{Betsumiya1}
Let $q$ be a power of an odd prime $p$ and $k$ a positive divisor of $q-1$. Then we have the following:\\
(i) $(a,b)_{e}\leq \Big \lceil \frac{k}{2}\Big \rceil$ for all $a,b$ with $0\leq a,b < e$ if $p>\frac{3k}{2}-1$;\\
(ii) $(a,b)_{e}\leq \Big \lceil \frac{k}{2}\Big \rceil-1$ for each $a$ with $0\leq a < e$ if $k$ is odd and  $p>\frac{3k}{2}$;\\
(iii) $(0,0)_{e}\leq \Big \lceil \frac{k}{2}\Big \rceil-1$ if $p>\frac{3k}{2}$;\\
(iv) $(0,0)_{e}=2$ if $p$ is sufficiently large compared to $k$ and $6|k$;\\
(v) $(0,0)_{e}=0$ if $p$ is sufficiently large compared to $k$ and $6\nmid k$.
\end{theorem}  
P. van Wamelen \cite{Wamelen1} has characterized the Jacobi sums of order $e$ corresponding to any given generator $\gamma$ of $\mathbb{F}_{q}^{*}$. So far this is the most satisfactory solution of the problem, as it takes up the case when $e$ is any integer $\geq 3$ and $q$ is any prime power for which $q \equiv 1 \pmod e$ and he proved:
\begin{theorem} \cite{Wamelen1}
Let $e\geq 3$, $p$ a prime, $q=p^{r}\equiv 1 \pmod e$. Let $q=ef+1$. Let $g_{m}=gcd(e,m)$, $g_{n}=gcd(e,n)$, $g=gcd(e,m+n)$, $g_{0}= gcd(g_{m},g_{n})$. Let $\epsilon_{g}(k)=1$ if $g|k$, and $0$ if $g\nmid k$. There is a unique polynomial $H\in \mathbb{Z}[x]$ such that $H(x)=a_{0}+a_{1}x+a_{2}x^{2}+\dots+a_{e-1}x^{e-1}$ and the coefficients satisfy the following three conditions:
\\
1. (a)
\begin{equation*}
\sum_{j=0}^{e-1}a_{j}^{2}=q+g_{0}ef^{2}-f(g_{m}+g_{n}+g).
\end{equation*}
(b) For $k=1,2,\dots,e-1$
\begin{equation*}
\sum_{j=0}^{e-1}a_{j}a_{j-k}=\epsilon_{g_{0}}(k)g_{0}ef^{2}-\epsilon_{g_{m}}(k)fg_{m}-\epsilon_{g_{n}}(k)fg_{n}-\epsilon_{g}(k)fg
\end{equation*}
where we consider the subscripts of the $a$'s modulo $e$. \\
2. 
\begin{equation*}
\sum_{k=0}^{e-1}ka_{k}\equiv \begin{cases}
0 \pmod e \hspace*{31mm} \text{if } e \text{ is odd},\\
(q-1)/2(g_{m}+g_{n}) \pmod e \hspace*{2mm} \text{if } e \text{ is  even}.
\end{cases}
\end{equation*}
3. For every $d$ dividing $e$ let $B_{d}\in \mathbb{Z}[x]$ be such that its reduction modulo $p$ is the minimal polynomial of $\gamma^{(q-1)/d}$ over $\mathbb{F}_{p}$ and $\prod_{k\in D_{d}}B_{d}(\zeta_{d}^{k})$ is not divisible by $p^{2}$ in $\mathbb{Z}[\zeta_{d}]$. Then $H(\zeta_{d})$ must satisfy the following conditions:\\
(a) if none of $m$, $n$ and $m+n$ are divisible by $d$,
\begin{equation*}
q|H(\zeta_{d})\prod_{k\in D_{d}}B_{d}(\zeta_{d}^{k^{-1}})^{(s_{d,q}(mk)+s_{d,q}(nk)-s_{d,q}(mk+nk))/(p-1)},
\end{equation*}
where $k^{-1}$ is taken modulo $d$.\\
(b) if $m \equiv -n \not \equiv 0 \pmod d$
\begin{equation*}
H(\zeta_{d})=-\chi_{d}^{m}(-1),
\end{equation*}
(c) if exactly one of $m$ and $n$ are divisible by $d$
\begin{equation*}
H(\zeta_{d})=-1,
\end{equation*}
(d) if both $m$ and $n$ are divisible by $d$
\begin{equation*}
H(\zeta_{d})=q-2.
\end{equation*}
If $H$ is the unique polynomial satisfying these three conditions, then
\begin{equation*}
H(\zeta_{d})=J(\chi^{m},\chi^{n}).
\end{equation*}
\end{theorem}
\section*{Acknowledgment}
\noindent
The authors acknowledge Central University of Jharkhand, Ranchi, Jharkhand for providing necessary and excellent facilities to carry out this research. 
M H Ahmed would like to thank K. Chakraborty and A. Hoque for their kind invitation and hospitality at Harish-Chandra Research Institute, where this paper 
was finalized.

\end{document}